\documentclass[conference,a4paper]{IEEEtran}

\usepackage[T1]{fontenc}
\usepackage[utf8]{inputenc}

\newtheorem{lem}{Lemma}

\newtheorem{theorem}{Theorem}

\newtheorem{problem}{Problem}
\newcommand{\be}{\begin{equation}}
\newcommand{\ee}{\end{equation}}
\newcommand{\ba}{\begin{array}}
\newcommand{\ea}{\end{array}}

\newcommand{\bea}{\begin{eqnarray*}}
\newcommand{\eea}{\end{eqnarray*}}
\newcommand{\bean}{\begin{eqnarray}}
\newcommand{\eean}{\end{eqnarray}}
 \def \rot{\mathrm{rot}}
 \def \div{\mathrm{div}}
\usepackage{cite}
\usepackage[cmex10]{amsmath}
\usepackage{amssymb,amsfonts}
\usepackage{array}
\usepackage[caption=false,font=footnotesize]{subfig}
\usepackage{dblfloatfix}

\usepackage{algorithmic}
\usepackage{graphicx}
\graphicspath{{graphics/}{moregraphics/}}
\DeclareGraphicsExtensions{.pdf,.PDF,.jpeg,.JPEG,.jpg,.JPEG,.png,.PNG}

\usepackage{booktabs}

\usepackage{siunitx}

\usepackage{textcomp}
\usepackage{xcolor}

\def\BibTeX{{\rm B\kern-.05em{\sc i\kern-.025em b}\kern-.08em
    T\kern-.1667em\lower.7ex\hbox{E}\kern-.125emX}}
\begin{document}

\title{High Order Impedance Boundary Condition for the Three-dimensional Scattering Problem in Electromagnetism}
\author{\IEEEauthorblockN{ OUESLATI Soumaya}
\IEEEauthorblockA{\textit{Department of Mathematics}\\
\textit{University CY Cergy Paris}\\
Cergy-Pontoise, France\\
soumaya.oueslati1@cyu.fr}

\and
\IEEEauthorblockN{ DAVEAU Christian}
\IEEEauthorblockA{\textit{Department of Mathematics}\\
\textit{University CY Cergy Paris}\\
Cergy-Pontoise, France\\
christian.daveau@cyu.fr}
\and
\IEEEauthorblockN{AUBAKIROV Abil}
\IEEEauthorblockA{\textit{Department of Mathematics}\\
\textit{University CY Cergy Paris}\\
Cergy-Pontoise, France\\
abil.aubakirov@cyu.fr}
}

\maketitle

\begin{abstract}

In this paper, we propose a variational formulation with the use of high order impedance boundary condition (HOIBC) to solve the scattering problem.  We study the existence and uniqueness  of the solution. Then, a discretization
of this formulation is done with Rao-Wilton-Glisson (RWG). We give validations of the HOIBC obtained with a 3D MoM code that show the improvement
in accuracy over the standard impedance boundary condition (SIBC) computations.

\end{abstract}

\begin{IEEEkeywords}
 boundary element method, scattering problem, Lagrange multipliers, high order impedance boundary condition.
\end{IEEEkeywords}

\section{Introduction}\label{intro}
We consider the scattering problem of electromagnetic waves $(\textbf{E}, \textbf{H})$ by a perfect conducting body with a complex coating. Let $\Omega$ be a  bounded domain with a Lipschitz-continuous boundary $\Gamma$. Let $n$ be the unit normal vector to $\Gamma$ directed to the exterior of $\Omega$. We define $\Omega_+$ as the space of radiating electric fields $\textbf{E}$ solutions of Maxwell equations in $\mathbb{R}^3\backslash\overline{\Omega_-}$ which have a tangential trace on $\Gamma$.

Waves propagate with constant wave number $k$ in the exterior unbounded domain $\Omega$. The electric field $\textbf{E}$ satisfies the harmonic Maxwell equation
$$\nabla \times(\nabla\times\textbf{E} )-k^2 \textbf{E}=0 $$
While the related magnetic field is given by
$$\textbf{H}=\frac{1}{ik}\nabla \times\textbf{E} $$
An electric field is said to be radiating if it satisfies the
Silver-M\"{u}ller radiation condition:
\begin{equation*}
 \lim_{r\rightarrow\infty} r(\textbf{E}\times\textbf{n}_r + \textbf{H}) = 0,
\end{equation*}
where $r = |\mathbf{x}|$ and $\mathbf{n}_r = \dfrac{\mathbf{x}}{|\mathbf{x}|}$, $\mathbf{x} \in \mathbb{R}^3$.\\
The coating is modeled by the following impedance boundary condition

\begin{equation}\label{eq1.2}
\textbf{E}_{t} - Z (n \times \textbf{H}) = 0 \,\,\,\,\,on\,\,\,\,\Gamma.
\end{equation}
Here,  $Z$ is the impedance operator that depends on incident angle, subscript tg denotes tangent component on the surface defined as:
$$\textbf{E}_{t} = n \times(\textbf{E} \times n).$$
Generally, we take the constant impedance operator, known as standard or Leontovitch impedance boundary conditions \cite{R-S, Senior_1995, Lange-1995}.
Our problem writes as follows:
\begin{equation}\label{problem_1}
\left\{
 \begin{array}{llll}
 \nabla \times(\nabla\times\textbf{E} )-k^2 \textbf{E}=0 \; \Omega_+, \\
\nabla \times(\nabla\times\textbf{H} )-k^2 \textbf{H}=0 \; \Omega_+, \\
   	   \textbf{E}_{t} - Z (n \times \textbf{H}) = 0 \,\,\,\,\,on\,\,\,\,\Gamma.\\
  \lim_{r\rightarrow\infty} r(\mathbf{E}\times\mathbf{n}_r + \mathbf{H}) = 0.
\end{array}
\right.
\end{equation}

This paper deals with the implementation of high order impedance boundary condition coupled by an integral equation.

\section{high order impedance boundary condition}

In this section, we propose a new impedance boundary condition with integral operators. In \cite{R-S} the authors approximated impedance condition as a ratio of polynomials of second order 
of the sine of the incidence angle in spectral domain. They obtained for the case of a rotational invariant coating a high order impedance boundary condition in spatial domain equations.
We transformed these equations with integral operators to get:
\begin{equation}\label{HOIBC3D}
(I + b_1 L_D - b_2 L_R)\textbf{E}_{t} = (a_0 I + a_1 L_D - a_2 L_R)(\mathbf{n} \times \textbf{H}).
\end{equation}
This approximation contains the  operators $L_D$ and $L_R$ are defined: for all vector function $\mathbf{A}$ sufficiently smooth, such that $\mathbf{A} . \mathbf{n}=0$ 
$$L_D (\mathbf{A}) =\nabla_{\Gamma}(\div_{\Gamma} \mathbf{A}), \,\,\,\,\,\,\,\, L_R(\mathbf{A})= \mathbf{\rot}_{\Gamma}(\rot_{\Gamma} \mathbf{A}).$$
In \cite{Abil}, different methods to calculate these coefficients $(a_0, a_j, b_j)$ for \ref{HOIBC3D}  are presented. In the following, we want to establish a new variational formulation, applying a boundary integral method for the problem in 3D case with HOIBC.
\section{Variational formulation with HOIBC}
We introduce current densities $\textbf{J} $ and $\textbf{M}$ on the boundary $\Gamma$ as follows: 
$$\mathbf{M} = \mathbf{E} \times \mathbf{n} \ \ , \ \mathbf{J} = \mathbf{n} \times \mathbf{H}.$$
We finally obtain the variational formulation
\begin{problem}\label{prob_3D}
 Find $U = (\mathbf{J}, \mathbf{M}, \tilde{\mathbf{J}}, \tilde{\mathbf{M}}) \in V = [H^{-1/2}(div, \Gamma) \cap L^2(\Gamma)]^4$ and $\lambda
 = (\boldsymbol{\lambda}_J, \boldsymbol{\lambda}_M) \in [H^{-1/2} (\Gamma)]^2$ such that 
 \begin{equation}\label{eq_prob_3D}
 \begin{cases}
 A(U, \Psi) + B^T(\lambda, \Psi) = F(\Psi) \\
   B(U, \lambda') = 0
 \end{cases}
 \end{equation}
 for all $\Psi = (\boldsymbol{\Psi}_J, \boldsymbol{\Psi}_M, \tilde{\boldsymbol{\Psi}}_J, \tilde{\boldsymbol{\Psi}}_M) \in V = [H^{-1/2}(div,
 \Gamma) \cap L^2(\Gamma)]^4$ and $\lambda' = (\boldsymbol{\lambda}'_J, \boldsymbol{\lambda}'_M) \in W = [H^{-1/2} (\Gamma)]^2$.
 \end{problem}
\noindent The bilinear forms are defined by:
 $$ B(U, \lambda') = \int_{\Gamma} \boldsymbol{\lambda}'_J \cdot (\tilde{\mathbf{J} } - \mathbf{n}\times\mathbf{J}) ds
 + \int_{\Gamma} \boldsymbol{\lambda}'_M \cdot (\tilde{\mathbf{M}} - \mathbf{n}\times \mathbf{M}) ds $$
 and
 \begin{equation}\label{oper_A_3d}
  A(U, \Psi) = iZ_0\iint_{\Gamma} kG\ (\textbf{J} \cdot \boldsymbol{\Psi}_J) - \frac{1}{k} G\ \div{\boldsymbol{\Psi}_J}\
  \div{\textbf{J}} ds ds'
 \end{equation}
 %
 %
 $$+ \frac{i}{Z_0} \iint_{\Gamma} kG\ (\boldsymbol{\Psi}_M \cdot \textbf{M}) - \frac{1}{k} G\ \div{\boldsymbol{\Psi}_M}\
 \div{\textbf{M}} dsds' $$
 $$+\iint_{\Gamma} \nabla'G \cdot (\boldsymbol{\Psi}_J \times \textbf{M}) dsds' - i \iint_{\Gamma} \nabla'G
 \cdot (\boldsymbol{\Psi}_M \times \textbf{J}) dsds'$$
 $$+ \frac{a_0}{2} \int_{\Gamma} \textbf{J} \cdot \boldsymbol{\Psi}_J ds + \frac{1}{2 a_0} \int_{\Gamma} \textbf{M}
 \cdot \boldsymbol{\Psi}_M ds $$
 $$ - \frac{a_1}{2} \int_{\Gamma} \div_{\Gamma}\textbf{J}\ \div_{\Gamma}\boldsymbol{\Psi}_J ds  - \frac{a_2}{2}
 \int_{\Gamma} \div_{\Gamma}\tilde{\textbf{J}}\ \div_{\Gamma} \tilde{\boldsymbol{\Psi}}_J ds$$
 $$ + \frac{b_1}{2} \int_{\Gamma} \div_{\Gamma}\tilde{\textbf{M}}\ \div_{\Gamma}\boldsymbol{\Psi}_J ds
 - \frac{b_2}{2} \int_{\Gamma} \div_{\Gamma}\textbf{M}\ \div_{\Gamma}\tilde{\boldsymbol{\Psi}}_J ds $$
 $$- \frac{b_1}{2 a_0} \int_{\Gamma} \div_{\Gamma}\tilde{\textbf{M}}\ \div_{\Gamma}\tilde{\boldsymbol{\Psi}}_M ds
 - \frac{b_2}{2 a_0} \int_{\Gamma} \div_{\Gamma}\textbf{M}\ \div_{\Gamma}\boldsymbol{\Psi}_M ds$$
 $$+ \frac{a_1}{2 a_0} \int_{\Gamma} \div_{\Gamma}\textbf{J}\ \div_{\Gamma}\tilde{\boldsymbol{\Psi}}_M ds
 - \frac{a_2}{2 a_0} \int_{\Gamma} \div_{\Gamma}\tilde{\textbf{J}}\ \div_{\Gamma}\boldsymbol{\Psi}_M ds, $$
 where $G(x,y)$ is the Green kernel.\\
 We use a theorem from \cite{JCN} to prove existence and uniqueness of a solution to the problem \ref{eq_prob_3D}. In the first time, we prove the following result.
 \begin{lem}\label{lemma_cont_A}
 The operator $A$ is continuous on $V\times V$ for all $\Psi \in V$ and we have that
  $$ | A(U, \Psi) |  \leq C \|U\|_V \| \Psi \|_V$$
 \end{lem}
\begin{lem}\label{lemma_coer_A}
The operator $A$ is coercive on $V$.
We have to show that there exist $\alpha > 0$ such that
 $$ \Re[A(U, U^*)] \geq \alpha \|U\|^2_V - C\|U\|^2_V, \ \forall U\in V.$$
 \end{lem}
 \begin{lem}\label{lemma_LBB}
The operator $B$ verifies the inequality
$$ \sup_{\|U\|_V=1}|B(U, \lambda)| \geq \beta \| \lambda \|_W,$$
$\forall \lambda \in W = [H^{-1/2} (\Gamma)]^3 \times [H^{-1/2} (\Gamma)]^3 $
where $U \in V = [H^{-1/2}(div, \Gamma) \cap L^2(\Gamma)]^4$ and $\beta > 0$.
\end{lem}
Therefore, we have the result.
\begin{theorem}
 The problem (\ref{eq_prob_3D}) admits a unique solution $U \in V = [H^{-1/2}(div, \Gamma) \cap L^2(\Gamma)]^4$ and $\lambda \in [H^{-1/2} (\Gamma)]^2$, if coefficients satisfy
 \begin{equation}\label{eq_ex&un_coef_3d}
  \Re(a_j) + \frac{|a_0||b_j + a^*_j/a^*_0|}{2} = 0 \ \ for\ \ j=1,2.
 \end{equation}
 \end{theorem}

\section{Discretization of the variational problem with HOIBC}
The first step consists to approach the surface of the obstacle by a surface $\Gamma_h$ composed of finite number of two dimensional elements.
These elements are triangular facets denoted by $T_ i$ for $i=1$ to $N_T$.
Let $\{ \mathbf{f}_i \}_{i=1,N_e}$ be a function of Rao-Wilton-Glisson functions. We decompose the electric and magnetic currents with thes basis functions.
For the  Lagrange multipliers, we use  basis functions proposed in \cite{Bendali_OCT_1999} whith
%
 is defined as follows: \\
\begin{equation}\label{function_g_k}
\mathbf{g}_n(x) =\left\{
\begin{array}{cc}
(1 - 2\omega^+_{i+2}(x) )(\nu_n \times \mathbf{n}^+),\ \ \ \ x \in T^+_n \\
(1 - 2\omega^-_{i+2}(x) )(\nu_n \times \mathbf{n}^-),\ \ \ \  x \in T^-_n
\end{array} \right.
\end{equation}
with $\nu_n$ is a direction vector of the edge $\mathbf{n}$ 
and $\{ \omega_i \}_{i = 1,3}$ are barycentric coordinates of $x$ relative to triangles $T^+_n$ or $T^+_n$.

Then, the discretized problem of  (\ref{eq_prob_3D}) is:
\begin{equation}\label{eq_sys_3D_h}
\begin{cases}
A^h(U_h, \Psi_h) + B^h(\lambda_h, \Psi_h) = \sum_{i=1}^{N_e} < \mathbf{E}^{inc}, \mathbf{f}_i > +\\ \sum_{i=1}^{N_e} < \mathbf{H}^{inc}, \mathbf{f}_i > \\
B^h(U_h, \lambda_h) = 0
\end{cases}
\end{equation}
where
$$ A^h(U_h, \Psi_h) = \sum^{N_e}_{i,j = 1} <Z_0 (B-S) \textbf{f}_j, \mathbf{f}_i > J_j + \sum^{N_e}_{i,j = 1} <Q \textbf{f}_j, \mathbf{f}_i> M_j  $$
$$ + Z^{-1}_0 \sum^{N_e}_{i,j = 1} < (B-S) \textbf{f}_j, \mathbf{f}_i > M_j- \sum^{N_e}_{i,j = 1} <Q\textbf{f}_j, \mathbf{f}_i > J_j + \frac{a_0}{2}$$
$$\sum^{N_e}_{i,j = 1} < \textbf{f}_j, \mathbf{f}_i> J_j + \frac{1}{2 a_0} \sum^{N_e}_{i,j = 1} < \textbf{f}_j, \mathbf{f}_i > M_j - \frac{a_1}{2}$$
$$\sum^{N_e}_{i,j = 1} < \mathbf{div}_{\Gamma}\textbf{f}_j, \mathbf{div}_{\Gamma}\mathbf{f}_i > J_j - \frac{a_2}{2} \sum^{N_e}_{i,j = 1} < \mathbf{div}_{\Gamma}\textbf{f}_j, \mathbf{div}_{\Gamma} \mathbf{f}_i > \tilde{J}_j  + \frac{b_1}{2}$$
$$\sum^{N_e}_{i,j = 1} < \mathbf{div}_{\Gamma}\textbf{f}_j, \mathbf{div}_{\Gamma}\mathbf{f}_i > \tilde{M}_j - \frac{b_2}{2} \sum^{N_e}_{i,j = 1} < \mathbf{div}_{\Gamma}\textbf{f}_j, \mathbf{div}_{\Gamma}\mathbf{f}_i > M_j- $$
$$ \frac{b_1}{2 a_0}\sum^{N_e}_{i,j = 1} < \mathbf{div}_{\Gamma}\textbf{f}_j, \mathbf{div}_{\Gamma}\mathbf{f}_i > \tilde{M}_j - \frac{b_2}{2 a_0} \sum^{N_e}_{i,j = 1} < \mathbf{div}_{\Gamma}\textbf{f}_j, \mathbf{div}_{\Gamma}\mathbf{f}_i > M_j$$
$$+ \frac{a_1}{2 a_0}\sum^{N_e}_{i,j = 1} < \mathbf{div}_{\Gamma}\textbf{f}_j, \mathbf{div}_{\Gamma}\mathbf{f}_i > J_j - \frac{a_2}{2 a_0} \sum^{N_e}_{i,j = 1} < \mathbf{div}_{\Gamma} \textbf{f}_j, \mathbf{div}_{\Gamma}\mathbf{f}_i > \tilde{J}_j$$
and
$$ B^h (U_h, \lambda'_h) = \sum^{N_e}_{i, k = 1} < \mathbf{g}_k, \mathbf{f}_i > \tilde{J}_i - \sum^{N_e}_{i,k = 1} < \mathbf{g}_k, \mathbf{n}\times\mathbf{f}_i > J_i$$
$$ + \sum^{N_e}_{i,k = 1} < \mathbf{g}_k, \mathbf{f}_i > \tilde{M}_i - \sum^{N_e}_{i,k = 1} < \mathbf{g}_k, \mathbf{n}\times\mathbf{f}_i > M_i $$
In the second step, we express the problem in matrix form
%
%
 \begin{equation*}
(B-S)_{i,j} = i\int \int_{\Gamma_h} k G(s,s') \mathbf{f}_j(s') \cdot \mathbf{f}_i(s)- \frac{1}{k} 
\end{equation*}
$$G(s,s') (\mathbf{div}_{\Gamma}\mathbf{f}_i) (\mathbf{div}_{\Gamma}' \mathbf{f}_j) ds ds' $$
$$Q_{i,j} = -i\int \int_{\Gamma_h} [\mathbf{f}_i(s) \times \mathbf{f}_j(s')] \cdot \nabla'_{\Gamma} G(s,s') ds ds' $$
$$I_{i,j} = \int_{\Gamma_h} \mathbf{f}_i \cdot \mathbf{f}_j ds ; \ \ \ \ D_{i,j} = \int_{\Gamma_h} (\mathbf{div}_{\Gamma}\mathbf{f}_j) (\mathbf{div}_{\Gamma}\mathbf{f}_i) ds $$
\begin{equation}\label{CH_CK}
C_{Hi,j} = \int_{\Gamma_h} \mathbf{g}_i \cdot \mathbf{f}_j ds ; \ \ \ \ C_{Ki,j} = \int_{\Gamma_h} \mathbf{g}_i \cdot (\mathbf{n} \times \mathbf{f}_j) ds,
\end{equation}

where $[C_H]$ is a non singular diagonal matrix. Therefore it is invertible. Then, we define $[A1]$ and $[A2]$ by
\begin{equation}\label{A1}
[A1] = [(B-S)] + \frac{a_0}{2}[I] - \frac{a_1}{2}[D],
\end{equation}
\begin{equation}\label{A2}
[A2] =[(B-S)] + \frac{1}{2 a_0}[I] - \frac{b_2}{2 a_0}[D].
\end{equation}
Using (\ref{CH_CK}), (\ref{A1}) and (\ref{A2}), we now present (\ref{eq_sys_3D_h}) in the following matrix form
\begin{equation}\label{system_ch_ck}
\left (
\begin{array}{cccccc}
[A1] & [Q] & 0 & \frac{b_1}{2}[D] & [C_K]^T & 0 \\
{[Q]}^T & [A2] & -\frac{a_2}{2 a_0}[D] & 0 & 0 & [C_K]^T \\
0 & -\frac{b_2}{2}[D] & -\frac{a_2}{2}[D] & 0 & [C_H]^T & 0 \\
\frac{a_1}{2 a_0}[D] & 0 & 0 & -\frac{b_1}{2 a_0}[D] & 0 & [C_H]^T \\
{[C_K]} & 0 & [C_H] & 0 & 0 & 0 \\
0 & [C_K] & 0 & [C_H] & 0 & 0
\end{array}
\right)
%
%
\end{equation}
\section{Elimination of auxiliary currents and the Lagrange multipliers}
Now, we focus ours attention on elimination auxiliary currents and the Lagrange multiplier. 
 Using the definition of the basis function (\ref{function_g_k}), we obtain that
$$\int_{\Gamma} \mathbf{g}_n(s) \cdot \tilde{\mathbf{J}} ds = \frac{|T^+_n| + |T^-_n|}{3} \tilde{J}_n .$$
On the other side, we have
$$\int_{\Gamma} \mathbf{g}_n(s) \cdot (\mathbf{n}\times\mathbf{J}) ds = \int_{T^+_n\cup T^-_n} \mathbf{g}_n(s) \cdot (\mathbf{n}\times\mathbf{J})  ds,$$
that is calculated with help of Gaussian quadrature.
And we get next equation
\begin{equation}\label{eq_chjt+ckj=0}
\frac{|T^+_n| + |T^-_n|}{3} \tilde{J}_n = \int_{T^+_n\cup T^-_n} \mathbf{g}_n(s) \cdot (\mathbf{n}\times\mathbf{J}) ds
\end{equation}
Thus, we obtain
\begin{equation}\label{eq_jt=-chinvckj}
\tilde{J}_n = \frac{3}{|T^+_n| + |T^-_n|} \int_{T^+_n\cup T^-_n} \mathbf{g}_n(s) \cdot (\mathbf{n}\times\mathbf{J}) ds
\end{equation}
on each edge of the mesh.
Consequently, we conclude that the auxiliary currents can be rewritten as
%
\begin{equation}\label{star}
\begin{array}{cc}
\overline{\tilde{J}} = -[C_H]^{-1} [C_K] \overline{J}\,\,\,\, \text{and} \,\,\,\,\overline{\tilde{M}} = -[C_H]^{-1} [C_K] \overline{M}
\end{array}
\end{equation}
%
So (\ref{system_ch_ck}) together with (\ref{star}) yields the following  expression of Lagrange multipliers in terms of $\overline{J}$ and $\overline{M}$
\begin{equation}\label{star_2}
\overline{\lambda_J} = \frac{b_2}{2} [C_H]^{-T} [D] \overline{M} - \frac{a_2}{2} [C_H]^{-T} [D] [C_H]^{-1} [C_K] \overline{J}
\end{equation}
\begin{equation}\label{star_3}
\overline{\lambda_M}= - \frac{a_2}{2 a_0} [C_H]^{-T} [D] \overline{J} - \frac{b_1}{2 a_0} [C_H]^{-T} [D] [C_H]^{-1} [C_K] \overline{M}
\end{equation}

%
%
%
%
However, the discretization  leads to dense linear systems.

The accurate computation of such matrix elements is an important. \cite{SICA} have been the evaluation of singular integrals.
\section{ Parallel BEM with H-matrices using MPI and OpenMP}

Because a dense coefficient HOIBC  matrix  described in the previous section, the total computational cost and the memory footprint of the method are very high. To remedy these problems, 
we developed a software framework that is based in parallel Hierarchical Matrices (H-matrices) with ACA (Adaptive Cross Approximation)  \cite{M_bebendorf}
\cite{SCB}.  
\begin{figure}[!h]
    \centering
        \includegraphics[scale=0.2]{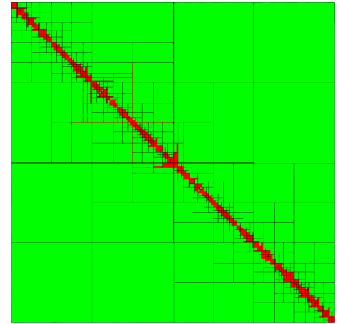}
        \caption{-matrix partition of a sphere mesh with 1434 edges.}
        \label{Hmatrix}
    \end{figure}
    
In Figure \ref{Hmatrix}, the matrix size is 2868x2868. Admissible blocks are shown in green, non-admissible blocks in red. \\
H-matrix computing can be conveniently coupled with parallel computing technologies to distribute the memory load, and to accelerate not only the matrix filling but also the classical matrix/vector product in order to obtain fast iterative solvers.
In next section, We present numerical results which show the effectiveness of both the parallelization and the H-matrix method. 
\section{Numerical results}

\begin{figure}[!h]
    \centering
        \includegraphics[scale=0.55]{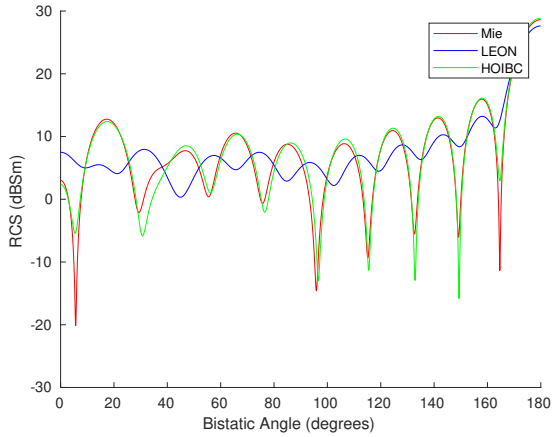}
        \caption{$\theta \theta$ component of the bistatic RCS for a coated conducting sphere.
        Exact serie solution and HOIBC solution.}
        \label{rot1}
    \end{figure}

We present now some numerical results obtained in 3D. 
As a first example, we consider the case of a coated conducting sphere with a conductor radius of $1.5 \lambda_0$ and coating thickness of $0.0075 \lambda_0$, with $\epsilon_r=5$ and $\mu_r=1.0$.
 Figure \ref{rot1} show the $\theta \theta$  components of the bistatic RCS for a plane wave incident from $\theta=0$. 
 Three solutions are included: the exact serie solution (MIE) and the solutions of the methods of moment studied below with SIBC and the HOIBC. 
 The Figure clearly shows the increased accuracy of the HOIBC solution relative to the SIBC solution. The SIBC gives only the average behavior of 
 the scattered field while the HOIBC accurately predicts the sidelobe behavior.
 
\begin{figure}[!h]
\centering
\includegraphics[scale=0.65]{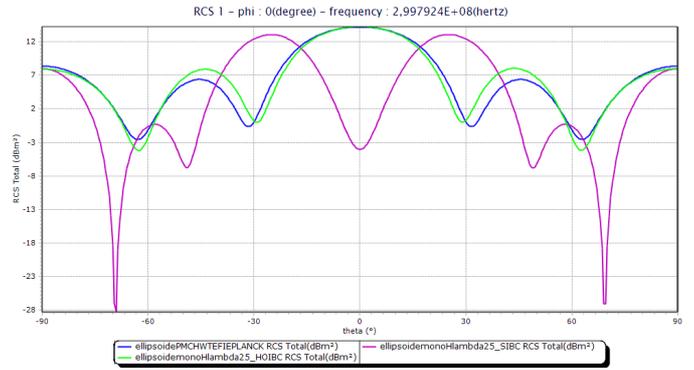}
\caption{$\theta \theta$ component of the monostatic RCS for a coated conducting spheroid, PMCHWT solution, SIBC solution and HOIBC solution.}
\label{rot3}
\end{figure}

\begin{figure}[!h]
\centering
        \includegraphics[scale=0.65]{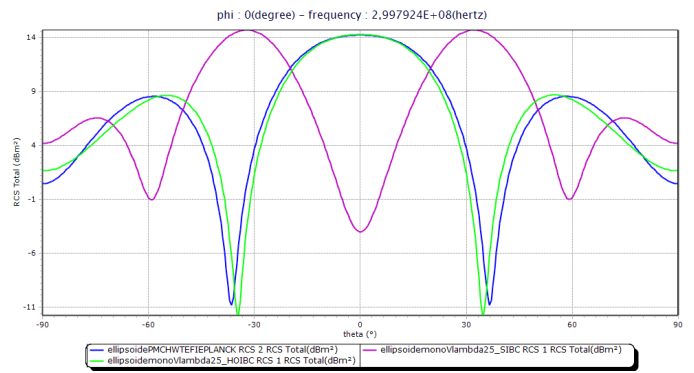}
        \caption{$\phi \phi$ component of the monostatic RCS for a coated conducting spheroid, PMCHWT solution, SIBC solution and HOIBC solution.}
        \label{rot4}
    \end{figure}

The second test is a coated conducting spheroid whose the radii are $0.5m$ and $1m$ with a coating thickness of $ 0.17 \lambda_0 $, with $ \epsilon_r = 5$ and $ \mu_r=1.0$.
Figures \ref{rot3} and \ref{rot4} show  the $\theta \theta$ and $\phi \phi$ components of the monostatic RCS. Three solutions are included : 
a method of moments solution called PMCHWT, the HOIBC solution and the SIBC solution. It does not exits exact solution for this case. 

We note that a slight difference between the PMCHWT solution and the HOIBC solution and the SIBC solution is very poor. In both cases, excellent results for HOIBC method are obtained 
for all angles of incidence.


\begin{figure}[!h]
    \centering
        \includegraphics[scale=0.35]{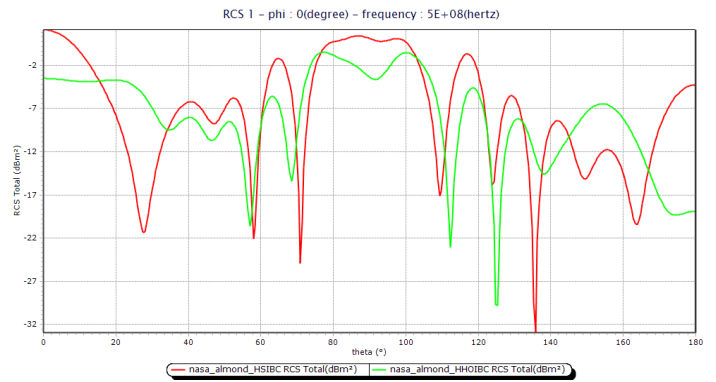}
        \caption{$\theta \theta$ component of the bistatic RCS for a coated conducting almond.
        Exact serie solution and HOIBC solution.}
        \label{almond}
    \end{figure}

As a final example a coated conducting almond is considered whose the radii are 0.25m and 2m with a coating thickness of $0.03$m, with $\epsilon_r=4$ and 
$\mu_r=1.0$. Figures \ref{almond} show  the $\theta \theta$ component of the bistatic RCS. Two solutions are included : the HOIBC solution and the SIBC solution. 

\section{conclusion}
The validations are shown an important accuracy improvement of the HOIBC model over the SIBC model in the case of materials of lower index and/or with small losses.

\end{document}